\documentclass[12pt, a4paper]{article}

 \usepackage{a4wide}
\usepackage{amsmath, amstext,amscd}

\usepackage{graphics}

\usepackage{amssymb}

\def\Xr{\mathcal{X}}
\def\JXr{\mathcal{JX}}
\def\Wr{\mathcal{W}}

\def\Vr{\widetilde{\mathcal{V}}}

\def\Z{\mathbb{Z}}

\def\O{\mathcal{O} }

\def\P{\mathbb{P}}
\def\Spec{\mathop{\rm Spec}\nolimits}
\def\div{\mathop{\rm div}\nolimits}
\def\Sym{\mathop{\rm Sym}\nolimits}
\def\Div{\mathop{\rm Div}\nolimits}
\def\Supp{\mathop{\rm Supp}\nolimits}
\def\Kum{\mathop{\rm Kum}\nolimits}
\def\GL{\mathop{\rm GL}\nolimits}

\makeatletter

\@addtoreset{equation}{section} \makeatother

\newtheorem{prop}{Proposition}[section]

\newtheorem{lem}[prop]{Lemma}
\newtheorem{theo}[prop]{Theorem}
\newtheorem{cor}[prop]{Corollary}

\def\Rq{\stepcounter{prop} \noindent \theprop{.} \emph{Remark.}  }

\def\dem{\noindent{\it Proof.} }

\def\tens{\otimes }

\def\Frac#1#2{{\displaystyle{{#1} \overwithdelims.. {#2}}}}

\begin{document}

\title{The action of the Frobenius map on rank $2$ vector bundles
over a
supersingular genus $2$ curve in characteristic $2$}
\author{Laurent Ducrohet}
\maketitle

\section{Introduction}

Let $X$ be a smooth proper genus $g$ curve over a field $k$ of
characteristic $p>0$. Denote by $M_X(r)$ the moduli space
of semistable vector bundles of rank $r$ and trivial determinant.
If $X_1$ is the curve $X\times_{k,\,\sigma} k$, where $\sigma : k
\to k$ is the Frobenius of $k$, the $k$-linear relative Frobenius
$F : X \to X_1$ induces by pull-back a (rational) map $V :
M_{X_1}(r) \dashrightarrow M_X(r)$. We are interested in
determining, e.g., the surjectivity and the degree of $V$, the
density of Frobenius-stable bundles, and the loci of
Frobenius-destabilized bundles, with the aim of studying the
behavior of the sequence $n \mapsto {F_{\rm abs}^{(n)}}^*E$ of
pull-back by $n$-fold iterated (absolute) Frobenius of a fixed
rank $r$ vector bundle $E$. Our interest in those questions comes
from a result (proved in [LS]) which claims that a semistable
rank $r$ vector bundle $E$ corresponds to an (irreducible)
continuous representation of the algebraic fundamental group
$\pi_1(X)$ in $\GL_r(\bar{k})$ (endowed with the discrete
topology) if and only if one can find an
integer $n>0$ such that ${F_{\rm abs}^{(n)}}^*E\cong E$.

For general $(g,\,r,\,p)$, not much seems to be known (see the
introductions of [LP1] and [LP2] for an overview of this subject).

When $g=2$, $r=2$, the isomorphism $D : M_X \to |2\Theta|\cong
\P^3$ (see [NR] for the complex case) remains valid for an
algebraically closed field of positive characteristic (see [LP1],
section 5 for a sketch of proof in the characteristic $2$ case) and
we have the commutative diagram
\begin{center}
\unitlength=0.6cm
\begin{picture}(10,4)
\put(1,3){$M_{X_1}(2)$}  \put(8,3){$M_{X}(2)$}
\put(1.2,0.5){$|2\Theta_1|$}  \put(8.2,0.5){$|2\Theta|$}
\put(3.4,3.2){\vector(1,0){4}} \put(3.4,0.7){\vector(1,0){4}}
\put(1.9,2.8){\vector(0,-1){1.5}}
\put(8.7,2.8){\vector(0,-1){1.5}} \footnotesize{\put(1.3,1.8){$D$}
\put(8.9,1.8){$D$} \put(5.2,3.4){$V$}
\put(5.2,0.9){$\widetilde{V}$}}
\end{picture} \end{center}
Furthermore, the semistable boundary of the moduli space $M_X(2)$
identifies via $D$ with the Kummer quartic surface $\Kum_X$,
which is canonically contained in the linear system $|2\Theta|$,
and $\widetilde{V}$ maps $\Kum_{X_1}$ to $\Kum_{X}$.
In [LP2], it is shown that $\widetilde{V}$ is given by degree $p$
polynomials
 and always has base-points.

When $p=2$ and $X$ is an ordinary curve, [LP1] determined the
quadric equations of $\widetilde{V}$ in terms of the generalized
theta constants of the curve $X$ and thus, they could answer the
above mentioned questions. In [LP2], they could give the equations
of $\widetilde{V}$ in case of a nonordinary curve $X$ with
Hasse-Witt invariant equal to $1$ by specializing a family $\Xr$
of genus $2$ curves parameterized by a discrete valuation ring
with ordinary generic fiber and special fiber isomorphic to $X$.
In particular, they determined the coefficients of the quadrics of
$\Vr_\eta$, which coincide with the Kummer quartic surface
coefficients, in terms of the coefficients of an affine equation
for a birational model of the ordinary curve $\Xr_\eta$.

In this paper, we complete the study of the
$(g,\,r,\,p)=(2,\,2,\,2)$ case by giving the equations of
$\widetilde{V}$ in case of a supersingular curve $X$ (Theorem
4.2). We adapt the strategy of [LP2]. Namely, we consider a family
$\Xr$ of genus $2$ curves parameterized by a discrete valuation
ring $R$ with ordinary generic fiber and special fiber isomorphic
to $X$. In order to find an $R$-basis of the free $R$-module
$\Wr:=H^0(\JXr,\,2\Theta)$ and to express the canonical theta
functions of order $2$ in that basis, we pull back $2\Theta$ on the
product $\Xr\times \Xr$ via the Abel-Jacobi map $\Xr \times \Xr
\to \JXr$ and make use of an explicit theorem of the square for
hyperelliptic Jacobians (see [AG]). As for the two other cases, we
can easily deduce a full description of the Verschiebung $V :
M_{X_1}(2) \dashrightarrow M_{X}(2)$ (Proposition
5.1).

I would like to thank C. Pauly for helpful discussions and Y.
Laszlo for having introduced me to this question, for his help and
encouragements.

\section{Deformation of genus $2$ curves}

\subsection{Specializing an ordinary curve}

Let $k$ be an algebraically closed field of characteristic $2$,
let $R$ be the discrete valuation ring $k[[s]]$, and let $K$ be its
fraction field. We consider a proper, smooth, and supersingular
curve $X$ of genus $2$ over $k$. By [L], we know that there exists
a unique $\mu \in k$ such that $X$ is birationally equivalent to
the plane curve given by the equation
\begin{eqnarray} y^2+y& = & x^5+\mu^2x^3\end{eqnarray} The
projection $(x,\,y) \mapsto x$ is the restriction of the ramified
double cover $\pi : X \to \P^1_k=|K_X|$ (where $K_X$ is the
canonical divisor of $X$) with a single Weierstrass point,
namely $\infty$.

Let us choose an element $\omega$ of $R-\{0,\,1\}$ and denote by
$\Xr$ the $R$-scheme defined by the two affine charts
\begin{eqnarray} y^2+(s^2x+1)(s^2\omega^2x+1)y=x^5+\mu^2x^3\end{eqnarray}
and
$$\widetilde{y}^2+\widetilde{x}(s^2+\widetilde{x})(s^2\omega^2+\widetilde{x})\widetilde{y}=\widetilde{x}+\mu^2\widetilde{x}^3$$
glued by the isomorphism given by $x\mapsto {\widetilde{x}}^{-1}$
and $y \mapsto \widetilde{y}(\widetilde{x}^{-3})$. The $R$-scheme
$\Xr$ is proper and smooth, and again, the projection $(x,\,y)
\mapsto x$ is the restriction of an $R$-morphism $\Xr \to \P^1_R$,
still denoted by $\pi$. It is convenient to introduce the notation
$$\left\{\begin{array}{l} q(x)=(s^2x+1)(s^2\omega^2x+1)=s^4\omega^2x^2+s^2(1+\omega^2)x+1\\
p(x)=x^5+\mu^2x^3\end{array}\right.$$

The special fiber $\Xr_0$ is isomorphic to $X$, the generic fiber
$\Xr_\eta$ is a proper and smooth ordinary curve of genus $2$ over
$K$, and $\pi_\eta : \Xr_\eta \to \P^1_K$ is a ramified double
cover with Weierstrass points $0_\eta$ (with coordinate $1/s^2$),
$1_\eta$ (with coordinate $1/s^2\omega^2$) and $\infty$. Thus, we
have $\JXr[2]_\eta \cong (\Z/2\Z)^2\times \mu^2_2/K$ and
$\JXr[2]_0$ is a self dual local-local group scheme over $k$ of
dimension $2$ and height $4$.

One can define (over $R$) the  Abel-Jacobi map
$$\begin{array}{cccl}
AJ : & \Xr \times \Xr & \to & \JXr\\
& (x,\,y) & \mapsto & \O(x+y)\tens K_{\Xr}^{-1}\end{array}$$ as
the compositum of the quotient map $\Xr \times \Xr \to
\Sym^2\Xr=\Div^2(\Xr)$ under the natural action of
$\mathfrak{S}_2$ by the natural map $\Div^2\Xr \to \JXr$.
It follows from the Riemann-Roch theorem that it is surjective, separable of
degree $2$. Denote by
$$[0]_\eta=AJ(1_\eta+\infty)\hskip8mm \ [1]_\eta=AJ(0_\eta+\infty)\hskip8mm
[\infty]_\eta=AJ(0_\eta+1_\eta)$$ the three nonzero elements of
$\JXr[2]^{\rm et}_\eta$ and notice that they specialize to $0$.

\subsection{Standard birational  model}

There exists a unique triplet $(a,\,b,\,c)$ of elements in the
algebraic closure $\bar{K}$ of $K$ such that $\Xr_\eta$ is
birationally equivalent to the plane curve given by the equation
\begin{eqnarray}
w^2+z(z+1)w= z(z+1)(az^3+(a+b)z^2+cz+c)
\end{eqnarray} and such
that, with these coordinates, $0_\eta=0$, $1_\eta=1$ and
$\infty=\infty$. These latter conditions amount to an ordering of
the Weierstrass points of $\Xr_\eta$, and this is equivalent to an
isomorphism $\mathbb{F}_2^2 \xrightarrow{\sim}
\JXr[2]^{\rm et}_\eta$, i.e., a level-$2$-structure on $X_\eta$ (see
[DO] for the characteristic $0$
case).

One can determine explicitly this coordinate change. As $x$ and
$z$ are rational coordinates on $\P^1_K$, the function $z=z(x)$
has to induce a $K$-isomorphism of $\P^1_K$, which fixes $\infty$,
sends $1/s^2$ to $0$ and $1/s^2\omega^2$ to $1$. Thus, we have
\begin{eqnarray}
z(x)=\frac{\omega^2}{1+\omega^2}(s^2x+1) & \Leftrightarrow &
x(z)=\frac{1+\omega^2}{s^2\omega^2}z+\frac{1}{s^2} \end{eqnarray}
In particular, one has
$q(x)=\Frac{(1+\omega^2)^2}{\omega^2}z(z+1)$ and replacing in
(2.2), one sees that the function $y(z,\,w)$ has to be of the form
\begin{eqnarray} y(z,\,w)=\frac{(1+\omega^2)^2}{\omega^2}[w+(\alpha
z^2+\beta z +\gamma)]\end{eqnarray}with $\alpha,\,\beta$ and
$\gamma$ in $\bar{K}$. Conversely, we have
\begin{eqnarray} w(x,\,y)  = \frac{\omega^2}{(1+\omega^2)^2}\left[y+
\alpha \,q(x) +(1+\omega^2)\left[(\alpha+\beta) s^2x+\left((\alpha
+\beta +\gamma)+\frac{\gamma}{\omega^2}\right)\right]\right]
\end{eqnarray}\\

\Rq At the Weierstrass point $0_\eta$, (2.2) gives
$y(0_\eta)^2=p\left(\Frac{1}{s^2}\right)$ and (2.3) gives
$w(0_\eta)=0$. Therefore, one has, using (2.5),
\begin{eqnarray}
\gamma^2=
\Frac{\omega^4}{(1+\omega^2)^4}\ p\left(\Frac{1}{s^2}\right)& =  &
\Frac{\omega^4}{s^{10}}\Frac{1+\mu
^2s^4}{(1+\omega^2)^4}\end{eqnarray}
 Similarly with the Weierstrass point $1_\eta$, one gets
\begin{eqnarray*}
(\alpha+\beta+\gamma)^2=\Frac{\omega^4}{(1+\omega^2)^4}\ p\left(\Frac{1}{s^2\omega^2}\right)&
=
&\Frac{1}{s^{10}\omega^{6}}\Frac{1+\mu^2s^4\omega^4}{(1+\omega^2)^4}\end{eqnarray*}
For further use, we observe
\begin{eqnarray} (\alpha+\beta)^2& = &
\Frac{1}{s^{10}\omega^6}\left[\Frac{1+\omega^{10}}{(1+\omega^2)^4}+
\mu^2s^4\omega^4\Frac{1+\omega^6}{(1+\omega^2)^4}\right]\\
(\alpha+\beta+\gamma)^2+\Frac{\gamma^2}{\omega^4} & =&
\Frac{1}{s^{10}\omega^6}\left[\Frac{1+\omega^{6}}{(1+\omega^2)^4}
+\mu^2s^4\omega^4\Frac{1}{(1+\omega^2)^3}\right]\hspace{0.3cm}
\end{eqnarray}
and
\begin{eqnarray}
(\alpha+\beta)^2+(1+\omega^2)^2\left((\alpha+\beta+\gamma)^2+\Frac{\gamma^2}{\omega^4}\right)&=&
\Frac{1}{s^{10}\omega^2}\Frac{1+\mu^2s^4\omega^2}{(1+\omega^2)^3}
\end{eqnarray}\\

Replacing $x$ and $y$ by their expressions (2.4) and (2.5), we
deduce from (2.2) the equation
\begin{eqnarray*}
w^2+z(z+1)w =\frac{\omega^4}{(1+\omega^2)^4}\left[
[(\alpha+\alpha^2)z^4+(\alpha+\beta)z^3+(\beta^2+\beta+\gamma)
z^2+\gamma z +\gamma^2]+p(x(z))\right]\end{eqnarray*} Identifying
the right-hand term of this equality and the right-hand term of
(2.3) (both are degree 5 polynomials), we get six equalities, the
first five of which are
\begin{eqnarray} \left\{\begin{array}{l}
a=\Frac{1+\omega^2}{s^{10}\omega^6}\\
b=\alpha^2+\alpha +\Frac{1}{s^{10}\omega^4}\\
a+b+c=\alpha+\beta+\Frac{\mu^2}{s^{6}\omega^2(1+\omega^2)}\\
0=\beta^2+\beta+\gamma +\Frac{\mu^2}{s^6(1+\omega^2)^2}\\
c=\gamma +\Frac{\omega^2(1+\mu^2s^4)}{s^{10}(1+\omega^2)^3}\\
\end{array}\right.
\end{eqnarray}
Equation (2.7) gives $\gamma$ and then, the last equality above
gives \begin{eqnarray}c =\Frac{1}{s^{10}}
\Frac{\omega^2}{(1+\omega)^6}(1+\nu_0)\text{ with
}\nu_0=\mu^2s^4+s^5(1+\mu s^2)(1+\omega^2).\end{eqnarray} Using
the expressions of $a$ and $c$, as well as the expression (2.9) of
$\alpha+\beta$, the third equality gives \begin{eqnarray}b
=\Frac{1}{s^{10}}\Frac{1}{\omega^6(1+\omega)^6}(1+\nu_1)\text{
with }\nu_1=\mu ^2s^4\omega^4+s^5\omega^3(1+\omega^2)(1+\mu
s^2\omega^2).\end{eqnarray}\\

\Rq Notice that the scalars $a,\,b$ and $c$ lie in fact in a
finite purely inseparable degree $2$ extension of the subfield of
$K$ generated by the coefficients in the equation (2.2). That is
why, for sake of simplicity, we have chosen Weierstrass points
$1/s^2$ and
$1/s^2\omega^2$ (which are squares in $K$) and the parameter $\mu^2$ in equation (2.1).\\

\Rq The second and the fourth equations of (2.11) show that the
coefficients $\alpha$ and $\beta$ lie in the finite extension
$K^{1/2}$ of $K$ and that they are not uniquely defined (but their
sum $\alpha+\beta$ is, and belongs to $K$). This indeterminacy
follows from the fact that the conditions we impose only determine
an automorphism of the projective line $\P^1_K$ and thus, an
automorphism of $\Xr_\eta$, defined up to the exchange of the two
sheets of the ramified cover $\pi : \Xr_\eta \to \P^1_K$, i.e., up to
  composition by the hyperelliptic involution of $\Xr_\eta$
(which is defined by the transformation $(z,\,w) \mapsto
(z,\,w+z(z+1)$)). More precisely, if $\alpha$ is a root of the
second equation of (2.11), $\alpha+1$ is the other one and,
changing $\alpha$ (resp. $\beta$) into $\alpha+1$ (resp.
$\beta+1$) in (2.5), we find that
\begin{eqnarray*} y'(z,\,w)& : = &
\frac{(1+\omega^2)^2}{\omega^2}[w+((\alpha+1)
z^2+(\beta+1) z +\gamma)]\\
 & =& \frac{(1+\omega^2)^2}{\omega^2}[(w+z(z+1))+(\alpha
z^2+\beta z +\gamma)]\\ & =& y(z,\,w+z(z+1)).\end{eqnarray*}

\section{The space $H^0(\JXr,\,\O(2\Theta))$ for a genus 2 curve $\Xr/R$}

Let $L$ be a field and $C$ a smooth proper curve of genus $2$ over
$L$. If $JC$ is its associated Jacobian variety, the Abel-Jacobi
map $AJ : C \times C \to JC$  is a surjective, separable morphism
of degree $2$ of $L$-varieties. The canonical divisor $K_C$ gives
a two-sheeted ramified cover $\pi : C \to \P^1_L$, and we denote
by $\bar{p}$ the involutive conjugate of any point $p$ in $C$.

It is classical that the natural map $\Sym^2C
=\Div^2(C) \to JC$ can be identified with the blowing up
$\widetilde{JC} \to JC$ of $JC$ at the origin: $\Div^2(C)
\to JC$ is a birational morphism of (nonsingular, projective)
surfaces and, using the Riemann-Roch theorem for $C$, we see that
$\bar{\Delta}$ is the only irreducible curve in $\Div^2(C)$
contracted to a point in $JC$, namely the origin. Using [H] (Chapter
V, Corollary 5.4.), we find that $\Div^2(C) \to JC$ is a
monoidal transformation, which is necessarily the monoidal
transformation $\widetilde{JC} \to JC$ with center the origin
since $\Div^2(C) \to JC$ factors through it (\emph{ibid.}
Proposition 5.3.).

 Suppose there exists an $L$-rational
ramification point $\infty$. Then, one can embed (noncanonically)
the curve $C$ in $JC$ by associating to each point $p\in C$ the
degree $0$ line bundle $\O(p-\infty)$. Denote by $\Theta$ the
corresponding divisor on $JC$. Its support is the set $\{j\in
JC\mid H^0(j\tens \O_C(\infty))\neq 0\}$ and it is well-known that
$\O(\Theta)$ is a principal polarization on the jacobian $JC$.
Consequently, we have $\dim H^0(JC,\,\O(2\Theta))=4$.

When considering the genus $2$ curve $\Xr \to \Spec R$ constructed
in the previous section, we can extend the construction and
therefore, $\Wr:=H^0(\JXr,\,\O(2\Theta))$ is a rank $4$ free
$R$-module.

\subsection{Canonical theta functions}

In the case of an ordinary genus $2$ curve over an algebraically
closed field of characteristic $2$, $\pi$ has three ramification
points. Upon composing with an automorphism of $\P^1_L$, one can
assume that they are $0,\,1$ and $\infty$, and that $C$ is
birational to a plane curve of the standard form (3), namely
$$w^2+z(z+1)w=z(z+1)(az^3+(a+b)z^2+cz+c)$$
with $a,\,b,\,c$ in $L$. Thus, the three nonzero $2$-torsion points
of $JC[2]$ are
$$[0]=AJ(1+\infty)\hskip8mm
[1]=AJ(0+\infty)\hskip8mm [\infty]=AJ(0+1)$$ Furthermore, $JC$ has a
canonical polarization $\Theta_B$ defined (see [R]) by means of
the canonical theta-characteristic $B \cong
\O_{C}(0+1+\infty)\tens K_{C}^{-1}$.

The $4$-dimensional $L$-vector space $W:=H^0(JC,\,\O(2\Theta_B))$
is the unique irreducible representation of weight $1$ of the
Heisenberg group $G(\O(2\Theta_B))$, obtained as a central
extension $$0 \to \mathbb{G}_m \to G(\O(2\Theta_B)) \to JC[2] \to
0$$ (\emph{cf.} [Mu2] for the general theory of Heisenberg groups
and [Sek] for the characteristic $2$ case).

Taking a nonzero section $\theta$ of $H^0(JC,\,\O(\Theta_B))$ and
setting
$$X_B=\theta^2\hskip8mm X_0=[0]^*X_B\hskip8mm X_1=[1]^*X_B\hskip8mm
X_\infty=[\infty]^*X_B$$ one obtains a basis (unique up to scalar
if we ask for these conditions) of $W$ (\emph{cf.} [LP1], section
2). As $\Supp(X_B)=\{j\in JC\mid H^0(j\tens B)\neq 0\}$, we have
$\Supp (X_\infty)=\{j\in JC \mid H^0(j\tens [\infty]\tens B)\neq
0\}$. But $[\infty]\tens B=\O_C(0+1)\tens \O_C(0+1+\infty) \tens
K_C^{-2}$ and using the fact that $0$ and $1$ are Weierstrass
points of $C$, one finds $[\infty]\tens B=\O_C(\infty)$. Thus,
$\Supp(X_\infty)=\Supp(\Theta)$ and finally
$2\Theta=[\infty]^*(2\Theta_B)$.

Following [LP2] (Lemma 3.3), one can express, using the
Abel-Jacobi map and an explicit theorem of the square ([AG]), the
rational functions
$\Frac{X_B}{X_\infty},\,\Frac{X_0}{X_\infty},\,\Frac{X_1}{X_\infty}
\in \bar{k}(JC)\subseteq \bar{k}(C\times C) $. We have the
following equalities in $\bar{k}(C\times C)$:
$$\begin{array}{l}F_B:=AJ^*\left(\Frac{X_B}{X_\infty}\right)=\alpha_B\Frac{(W_1+W_2)^2}{P(z_1,\,z_2)}
\text{ where }\left\{\begin{array}{l}W_i=\Frac{w_i}{z_i(z_i+1)}\\
P(z_1,\,z_2)=\Frac{(z_1+z_2)^2}{z_1z_2(z_1+1)(z_2+1)}\end{array}\right.\\
F_0:=AJ^*\left(\Frac{X_0}{X_\infty}\right)=\alpha_0z_1z_2,
\hspace{0.5cm}
F_1:=AJ^*\left(\Frac{X_1}{X_\infty}\right)=\alpha_1(z_1+1)(z_2+1)\end{array}$$
with nonzero scalars $\alpha_B,\,\alpha_0$ and $\alpha_1$
explicitly determined in terms of the scalars $a,\,b,\,c$
appearing in the equation of the standard birational model of $C$,
namely
$$\alpha_B=\Frac{1}{\sqrt{bc}}\hspace{8mm} \alpha_0=\sqrt{\Frac{a}{c}} \hspace{8mm}
\alpha_1=\sqrt{\Frac{a}{b}}$$

\subsection{Case of the genus 2 curve $\Xr \to \Spec R$}

The latter results only apply over the generic fiber $\Xr_\eta$
but we can use the coordinates change formulae (2.4) and (2.6) to
express the rational functions $F_B,\,F_0$ and $F_1$ in terms of
the $R$-coordinates $(x_i,\,y_i)$ instead of  the $K$-coordinates
$(z_i,\,w_i)$. First, using the expressions (2.11), (2.12) and
(2.13) of coefficients $a,\,b$ and $c$, and putting $\tau_0$ and
$\tau_1$ for $\sqrt{\Frac{1}{1+\nu_0}}$ and
$\sqrt{\Frac{1}{1+\nu_1}}$ respectively (which belong to
$R[\sqrt{s}]$), we have
\begin{eqnarray}\sqrt{\Frac{a}{c}}=\Frac{1+\omega^4}{\omega^4}\tau_0\hspace{0.8cm}
\sqrt{\Frac{a}{b}}=(1+\omega^4)\tau_1\hspace{0.8cm}
\sqrt{\Frac{1}{bc}}=s^{10}\omega^2(1+\omega^2)\tau_0\tau_1\end{eqnarray}
Then, straightforward computations give
\begin{eqnarray}
F_0=\tau_0(1+s^2(x_1+x_2)+s^4x_1x_2)\hspace{8mm}
F_1=\tau_1(1+s^2\omega^2(x_1+x_2)+s^4\omega^4x_1x_2)
\end{eqnarray}
and
\begin{eqnarray*} F_B& =&\tau_0\tau_1\left[
s^6\omega^2(1+\omega^2)\Frac{1}{q(x_1)q(x_2)}\left(\Frac{q(x_2)y_1+q(x_1)y_2}{x_1+x_2}\right)^2\right.\\
& & \hspace{1cm}+
s^{18}\omega^6(1+\omega^2)^3(\alpha+\beta)^2\Frac{(x_1x_2)^2}{q(x_1)q(x_2)}\\
& & \hspace{1.2cm}
+\left(s^{14}\omega^6(1+\omega^2)^3\left((\alpha +\beta
+\gamma)^2+\Frac{\gamma^2}{\omega^4}\right)\right)\Frac{(x_1+x_2)^2}{q(x_1)q(x_2)}\\
& &
\left.\hspace{1.4cm}+s^{10}\omega^2(1+\omega^2)^3\left((\alpha+\beta)^2+(1+\omega^2)^2\left((\alpha
+\beta
+\gamma)^2+\Frac{\gamma^2}{\omega^4}\right)\right)\Frac{1}{q(x_1)q(x_2)}\right]
\end{eqnarray*}
We compute the coefficients of the last three terms, using
equalities (2.8), (2.9) and (2.10) and finally obtain
\begin{eqnarray}\begin{array}{rcl} F_B& =&\tau_0\tau_1\left[
s^6\omega^2(1+\omega^2)\Frac{1}{q(x_1)q(x_2)}\left(\Frac{q(x_2)y_1+q(x_1)y_2}{x_1+x_2}\right)^2\right.\\
& & \hspace{1cm}+ s^{8}\left(\Frac{1+\omega^{10}}{1+\omega^2}+
\mu^2s^4\omega^4\Frac{1+\omega^6}{1+\omega^2}\right)\Frac{(x_1x_2)^2}{q(x_1)q(x_2)}\\
& & \left.\hspace{1.2cm}
+s^{4}\left(\Frac{1+\omega^{6}}{1+\omega^2}+
\mu^2s^4\omega^4\right)\Frac{(x_1+x_2)^2}{q(x_1)q(x_2)}+(1+\mu^2s^4\omega^2)\Frac{1}{q(x_1)q(x_2)}\right]
\end{array}\end{eqnarray}

\Rq Notice that, as the product $q(x_1)q(x_2)$ specializes to $1$
over the special fiber, the three rational functions
$\Frac{X_0}{X_\infty},\,\Frac{X_1}{X_\infty}$ and
$\Frac{X_B}{X_\infty}$ specialize to $1$ over the special fiber.
This corresponds to the specialization of the $2$-torsion points
$[0]_\eta,\,[1]_{\eta}$ and
$[\infty]_\eta$ to $0$.

\subsection{Finding an $R$-basis}

Let us consider again a  smooth   proper curve $C$ of genus $2$ over  an arbitrary field
$L$. Denote by $\imath$ the hyperelliptic involution that permutes
the two sheets of the ramified cover $\pi : C \to \P_L^1$ given by
the canonical map and assume that $\pi$ has a $L$-rational
ramification point $\infty$. In the sequel, we will use the
Abel-Jacobi map $AJ:  C \times C \to JC$ and the divisor $\Theta$
on $JC$ defined as the image of $C$ in $JC$ via the embedding $p
\mapsto \O_C(p-\infty)$.

\begin{lem} The pull-back $AJ^*(\Theta)$ is the divisor
$(\{\infty\} \times C+C\times \{\infty\})+\bar{\Delta}$ on
$C\times C$, where $\bar{\Delta}=(Id\times \imath)^*\Delta$ and
$\Delta$ is the diagonal in the product $C\times C$.\end{lem}

 \dem
Let $M$ (resp. $N$) be the prime divisor of $\Sym^2 C$
with support the set $\{(p+\infty),\,p\in C\}$ (resp. the set
$\{(p+\bar{p}),\,p\in C\}$). From [AG], we know that
$b^*(\Theta)=M+nN$, where $n$ is a nonnegative integer, and that
$\sigma^*(M)=C \times \{\infty\} + \{\infty\}\times C$.
Furthermore, as $\sigma^{-1}(N)=\bar{\Delta}$,
$\sigma^*(N)=k\bar{\Delta}$, where $k$ is an nonnegative integer
satisfying the equation $\deg(\sigma)(N)^2=k^2(\bar{\Delta})^2$.

On the one hand, since $N$ is the exceptional curve of the blowing
up $b : \Sym^2C \to JC$, we have $(N)^2=-1$. On the other
hand, the self-intersection number $(\bar{\Delta})^2$ is coincides
with $\deg_{\bar{\Delta}}(\O(\bar{\Delta})\tens
\O_{\bar{\Delta}})$. But $\O(\bar{\Delta})\tens
\O_{\bar{\Delta}}=\mathcal{N}_{\bar{\Delta}/C\times C}$ and taking
determinants in the following short exact sequence of
$\O_{\bar{\Delta}}$-coherent sheaves
$$0 \to \mathcal{N}_{\bar{\Delta}/C\times C}^{-1} \to \Omega_{C\times
C}^1 \tens \O_{\bar{\Delta}} \to \omega_{\bar{\Delta}} \to 0$$ we
obtain $\omega_{C}\tens \imath^*\omega_{C} \cong \omega_{C}^2
\cong \mathcal{N}_{\bar{\Delta}/C\times C}^{-1} \tens
\omega_{\bar{\Delta}}$ hence
\begin{eqnarray} \omega_{\bar{\Delta}}^{-1}\cong
\mathcal{N}_{\bar{\Delta}/C\times C}\end{eqnarray} Therefore,
$(\bar{\Delta})^2=2-2g_C=-2$, hence $k=1$ and $AJ^*(\Theta)=C
\times \{\infty\} + \{\infty\}\times C+n\bar{\Delta}$. One can
compute self-intersection again to determine $n$. We write
$$\begin{array}{rl}
\deg(AJ)(\Theta)^2=& (C \times \{\infty\})^2 +  (
\{\infty\}\times C)^2+n^2(\bar{\Delta})^2\\
 & + 2[(C \times \{\infty\}).(\{\infty\}\times C)+n((C \times
\{\infty\}).\bar{\Delta}+(\{\infty\}\times
C).\bar{\Delta})]\end{array}$$ It is clear that $(C \times
\{\infty\})^2$ (resp. $(\{\infty\}\times C)^2$) equals to zero for
$C \times \{\infty\}$ (resp. $\{\infty\}\times C$) being
algebraically equivalent to a divisor that does not meet $C \times
\{\infty\}$ (resp. $\{\infty\}\times C$). It is clear as well that
the intersection products $(C \times \{\infty\}).(\{\infty\}\times
C)$, $(C \times \{\infty\}).\bar{\Delta}$ et $(\{\infty\}\times
C).\bar{\Delta}$ equal to $1$ since the two divisors in each pairs
obviously meet transversally in a single point, namely
$(\infty,\,\infty)$. Furthermore, using the Riemann-Roch theorem
for an abelian variety of dimension $2$ and a principal (hence
ample)
divisor ([Mu1]), we have $(\Theta)^2=2$. Replacing, we obtain $n=1$. $\square$ \\

Denote by $p_i$ the canonical projection $C\times C \to C$ on the
$i$-th factor ($i=1,\,2$). Because the canonical divisor of $C$ is
$2\infty$ and because $\Omega^1_{C\times C}\cong p_1^*(\omega_{C})
\oplus p_2^*(\omega_{C})$, the canonical divisor $K_{C\times C}$
equals $2(C \times \{\infty\}) + 2(\{\infty\} \times C)$.
Therefore, the previous lemma allows us to see the $L$-vector
space $H^0(JC,\,\O(2\Theta))$ as the linear subspace of
$H^0(C\times C,\,\O(K_{C\times C}+2\bar{\Delta}))$ consisting of
 symmetric sections (under the action of
$\mathfrak{S}_2$)
that take constant value along $\bar{\Delta}$.

\begin{lem} The natural inclusion $H^0(C\times C,\,K_{C\times C})\hookrightarrow H^0(C\times
C,\,\O(K_{C\times C}+2\bar{\Delta}))$ induces three linearly
independent sections $1$, $x_1+x_2$ and $x_1x_2$ of
$H^0(JC,\,\O(2\Theta))$.\end{lem} \dem Consider the following
short exact sequence of $\O_{C \times C}$-coherent sheaves
\begin{eqnarray}0 \to \O(-\bar{\Delta}) \to \O_{C\times C} \to
\O_{\bar{\Delta}} \to 0\end{eqnarray} Tensoring with
$\O(K_{C\times C}+2\bar{\Delta})$ and taking cohomology groups,
one gets the following exact sequence of $L$-vector spaces
$$0 \to H^0(\O(K_{C\times
C}+\bar{\Delta})) \to H^0(\O(K_{C\times C}+2\bar{\Delta})) \to
H^0(\O(K_{C\times C}+2\bar{\Delta})\tens \O_{\bar{\Delta}})$$ But
$\O(K_{C\times C}+2\bar{\Delta})\tens \O_{\bar{\Delta}}\cong
\left(\O(K_{C\times C})\tens \O(\bar{\Delta})\tens
\O_{\bar{\Delta}}\right)\tens_{\O_{\bar{\Delta}}}\left(\O(\bar{\Delta})\tens
\O_{\bar{\Delta}}\right)$. The first term of this tensor product
is isomorphic to $\omega_{\bar{\Delta}}$ and the second to $
\mathcal{N}_{\bar{\Delta}/C\times C}$. Thus, using the isomorphism
(3.4), we obtain the structural sheaf $\O_{\bar{\Delta}}$ and we
have an exact sequence
$$0 \to H^0(\O(K_{C\times
C}+\bar{\Delta})) \to H^0(\O(K_{C\times C}+2\bar{\Delta})) \to L
$$ Tensoring (3.5) with $\O(K_{C\times C}+\bar{\Delta})$ and
taking cohomology groups, one gets an exact sequence of $L$-vector
spaces
$$0 \to H^0(\O(K_{C\times C})) \to H^0(\O(K_{C\times
C}+\bar{\Delta})) \to H^0(\O(K_{C\times C}+\bar{\Delta})\otimes
\O_{\bar{\Delta}}) \xrightarrow{\delta} H^1(\O(K_{C\times C}))$$
with $\O(K_{C\times C}+\bar{\Delta})\otimes \O_{\bar{\Delta}}\cong
\omega_{\bar{\Delta}}$. By Serre duality, the dual exact sequence
corresponds to the long exact sequence of cohomology groups
associated to the short exact sequence (3.5). Therefore, using
Serre duality again and K\"unneth isomorphism, one has the
following commutative diagram of morphisms of $L$-vector spaces

\begin{center}
\unitlength=0.6cm
\begin{picture}(12,6)
\put(0.6,0){\put(0,5.5){$H^1(\O(K_{C\times C}))^\vee$}
\put(9,5.5){$H^0(\omega_{\bar{\Delta}})^\vee$}
\put(4.8,5.6){\vector(1,0){3.4}}
 \put(10,3.5){\vector(0,1){1.6}} \put(10.2,4){$\wr$}
\put(2.2,3.5){\vector(0,1){1.6}} \put(2.4,4){$\wr$}
\put(0.6,2.8){$H^1(\O_{C\times C})$}
\put(9,2.8){$H^1(\O_{\bar{\Delta}})$}
\put(4.8,2.9){\vector(1,0){3.4}}
 \put(5.5,0.1){\vector(1,0){2.8}}
\put(10,0.7){\vector(0,1){1.6}} \put(10.2,1.2){$\wr$}
\put(2.2,0.7){\vector(0,1){1.6}} \put(2.4,1.2){$\wr$}
\put(9,0){$H^1(\O_{C})$} {\small \put(5.6,0.3){$Id+H^1(\imath)$}
\put(6,5.8){$\delta^\vee$}} }

\put(0,0){$H^1(\O_{C}) \oplus H^1(\O_{C})$}

\end{picture}
\end{center}
As the bottom horizontal arrow is surjective, $\delta^\vee$ is
surjective and $\delta$ is injective. Thus, the morphism
$H^0(\O(K_{C\times C})) \to H^0(\O(K_{C\times C}+\bar{\Delta}))$
is an isomorphism and using the K\"unneth isomorphism again, one
gets the following exact sequence of $L$-vector spaces
\begin{eqnarray*} 0 \to H^0(\omega_{C}) \tens H^0(\omega_{C})\to
H^0(\O(K_{C\times C}+2\bar{\Delta})) \to L
\end{eqnarray*}
Denote by $x$ the rational coordinate function of $\P^1_L$ with
pole at $\infty$. Considering $H^0(\omega_{C})$ as an $L$-vector
subspace of the function field $\bar{k}(C)$ by means of the
differential $dx$, it has a basis $\{1,\,x\}$. Thus, considering
$H^0(\omega_{C}) \tens H^0(\omega_{C})$ as an $L$-vector subspace
of the function field $\bar{k}(C\times C)$, we find that
$\left(H^0(\omega_{C}) \tens
H^0(\omega_{C})\right)^{\mathfrak{S}_2}$ has a basis
$\{1,\,x_1+x_2,\,x_1x_2\}$. Note that the three corresponding
sections of $H^0(\O(K_{C\times
C}+2\bar{\Delta}))^{\mathfrak{S}_2}$ are in the kernel of the
evaluation $H^0(\O(K_{C\times C}+2\bar{\Delta})) \to L\cong
H^0(\O_{\bar{\Delta}})$ so they have constant value (equal to
zero) along $\bar{\Delta}$, hence define sections of $H^0(JC,\,\O(2\Theta))$. $\square$\\

\Rq Consider the exact sequence (3.5), tensored by $\O(K_{C\times
C} + 2\bar{\Delta})$ and localized at the generic point
$\bar{\Delta}$ of the irreducible subscheme $\bar{\Delta}$ of
$C\times C$. If $A$ denotes the discrete valuation ring
$\O_{C\times C,\,\bar{\Delta}}$, if $t$ is a element of $A$ that
generates the maximal ideal, we find the exact sequence $$0 \to
t^{-1}.A \to t^{-2}.A \to \bar{k}(\bar{\Delta}) \to 0$$ Now, given
$g \in H^0(\O(K_{C\times C}+2\bar{\Delta}))$, we can look at it as
an element of $t^{-2}.A$ and write it as
$a_{-2}t^{-2}+a_{-1}t^{-1}+\sum_{n\geq 0}a_nt^n$, where the $a_i$
 are \emph{a priori} elements of the residue field
$\bar{k}(\bar{\Delta})$. Its class in $\bar{k}(\bar{\Delta})$ is
$a_{-2}$, which lies in fact in $L$ since the following diagram
(where the vertical arrows are localization) commutes

\begin{center}
\unitlength=0.6cm
\begin{picture}(12,4)
\put(0,2.5){$H^0(\O(K_{C\times C}+2\bar{\Delta}))$}
\put(3.2,2.3){\vector(0,-1){1.4}} \put(2.2,0.2){$t^{-2}.A$}
\put(1.2,0){\put(9,2.5){$L$}
\put(8.7,0.2){$\bar{k}(\bar{\Delta})$}
\put(3.5,0.4){\vector(1,0){4.5}} \put(5.5,2.7){\vector(1,0){2.5}}
\put(9.2,2.3){\vector(0,-1){1.4}}}
\end{picture}
\end{center}
Therefore, the morphism $H^0(\O(K_{C\times C}+2\bar{\Delta})) \to
L$ amounts, in
some sense, to the computation of a residue.\\

\Rq As $H^0(JC,\,\O(2\Theta))$ is of dimension $4$, the morphism
$H^0(\O(K_{C\times C}+2\bar{\Delta})) \to L$ cannot be zero and
there must exist a symmetric rational function $f$ on $C\times C$
such that
\begin{enumerate}
\item[(1)] $\div(f)+K_{C\times C}+2\bar{\Delta}$ is
effective;
 \item[(2)] $f$ has a pole of order $2$ along
$\bar{\Delta}$.\end{enumerate} Thus, we will have a basis
$\{1,\,x_1+x_2,\,x_1x_2,\,f\}$ of the vector space
$H^0(JC,\,\O(2\Theta))$.\\

\Rq These two lemmas extend to the relative case of the genus $2$
curve $\Xr$ over $\Spec R$, and we find three linearly
independent sections $1$, $x_1+x_2$ and $x_1x_2$ of $\Wr$.\\

We now restrict to the case of the genus $2$ curve $\Xr \to \Spec
R$. Note that $F_B$ certainly satisfies condition (1) in the
remark above: over the generic point, it is the pull-back by $AJ$
of the rational function $\Frac{X_B}{X_\infty}$, where $X_B$ and
$X_\infty$ are two sections of the line bundle $\O(2\Theta)$,
which corresponds, by Lemma 3.2, to the divisor $K_{\Xr_\eta\times
\Xr_\eta}+2\bar{\Delta}$. Furthermore, the origin of $\JXr_\eta$
does not belong to $\Supp(\Theta_B)$: if it did, we would
have $h^0(\Xr_\eta, B)\geq 1$ and there would exist a point $p$ in
$\Xr_\eta$ such that the divisors $0_\eta+1_\eta$ and $p+\infty$
were linearly equivalent, which is impossible since the three
Weierstrass points of $\Xr_\eta$ are pairwise different. On the
other hand, $H^0(\Xr_\eta, \O(\infty))\geq 1$ so
$\Frac{X_B}{X_\infty}$ has a pole at the origin, and by pull-back,
$F_B$ satisfies condition (2) over the generic fiber.

Of course, using Lemma 3.3 and Remark 3.5, any linear combination
$$\mu\left(\Frac{1}{\tau_0\tau_1}F_B+(\lambda+\lambda_{\Sigma} (x_1+x_2)+ \lambda_{\Pi}
(x_1x_2))\right)$$ with $\mu,\,\lambda,\,\lambda_{\Sigma}$ and $
\lambda_{\Pi}$ in $K$ such that it belongs to the ring
$$A=R[x_1,\,x_2,\,y_1,\,y_2,\,(x_1+x_2)^{-1},\,(q(x_1)q(x_2))^{-1}]$$
will fulfill condition (1) over both generic and special fibers
and condition (2) over the generic fiber.

Using (3.3), recall that
\begin{eqnarray*}\begin{array}{rcl} \Frac{1}{\tau_0\tau_1}F_B& =&
s^6\omega^2(1+\omega^2)\Frac{1}{q(x_1)q(x_2)}\left(\Frac{q(x_2)y_1+q(x_1)y_2}{x_1+x_2}\right)^2\\
& & \hspace{0.2cm}+ s^{8}\left(\Frac{1+\omega^{10}}{1+\omega^2}+
\mu^2s^4\omega^4\Frac{1+\omega^6}{1+\omega^2}\right)\Frac{(x_1x_2)^2}{q(x_1)q(x_2)}\\
& & \hspace{0.4cm} +s^{4}\left(\Frac{1+\omega^{6}}{1+\omega^2}+
\mu^2s^4\omega^4\right)\Frac{(x_1+x_2)^2}{q(x_1)q(x_2)}+(1+\mu^2s^4\omega^2)\Frac{1}{q(x_1)q(x_2)}
\end{array}\end{eqnarray*}
As
$\Frac{1}{q(x_1)q(x_2)}\left(\Frac{q(x_2)y_1+q(x_1)y_2}{x_1+x_2}\right)^2$
obviously has a pole of order $2$ along $\bar{\Delta}$ over both
generic and special fibers, we look for
$\lambda,\,\lambda_{\Sigma}$ and $ \lambda_{\Pi}$ in $K$ such that
we can take $\mu=\Frac{1}{s^6\omega^2(1+\omega^2)}$. Looking at
the class of $\Frac{1}{\tau_0\tau_1}F_B$ in $A/s^6$, we find
\begin{eqnarray*} \Frac{1}{q(x_1)q(x_2)}\left[
s^{4}\Frac{1+\omega^{6}}{1+\omega^2}(x_1+x_2)^2+(1+\mu^2s^4\omega^2)\right]
\end{eqnarray*}
On the other hand, $(\lambda+\lambda_{\Sigma} (x_1+x_2)+
\lambda_{\Pi} (x_1x_2))q(x_1)q(x_2)$ equals
$$\lambda+(\lambda_\Sigma+\lambda
s^2(1+\omega^2))(x_1+x_2)+(\lambda_\Pi+\lambda
s^4(1+\omega^4))(x_1x_2)+(\lambda s^4\omega^2+\lambda_\Sigma
s^2(1+\omega^2))(x_1+x_2)^2$$ in $A/s^6$. Thus, we see that
$$\Frac{1}{\tau_0\tau_1}F_B+(1+\mu^2s^4\omega^2)+s^2(1+\omega^2) (x_1+x_2)+
s^4(1+\omega^4) (x_1x_2)$$ belongs to $s^6A$. A few more
calculations show that it is in fact in $
s^6\omega^2(1+\omega^2)A$ and we set
\begin{eqnarray*}f=\Frac{1}{s^6\omega^2(1+\omega^2)}\left[\Frac{1}{\tau_0\tau_1}F_B+(1+\mu^2s^4\omega^2)+s^2(1+\omega^2)
(x_1+x_2)+ s^4(1+\omega^4) (x_1x_2)\right]\end{eqnarray*}
Therefore,
\begin{eqnarray}
F_B=\tau_0\tau_1((1+\mu^2s^4\omega^2)+s^2(1+\omega^2) (x_1+x_2)+
s^4(1+\omega^4) (x_1x_2)+s^6\omega^2(1+\omega^2)f)
\end{eqnarray}

\begin{prop} The projective space $\P \Wr \cong \P^3_R$ has
homogeneous coordinates $\{z_1,\,z_2,\,z_3,\,z_\infty\}$ in terms
of which, over the generic point $\Spec K \to \Spec
R$, the canonical theta coordinates
$\{x_B,\,x_0,\,x_1,\,x_\infty\}$ are given, up to scalar, by the
formulae
$$\left\{\begin{array}{l} x_B=z_1\\
x_0=\tau_1(z_1+s^2\omega^2z_2+s^4\omega^4z_3)\\
x_1=\tau_0(z_1+s^2z_2+s^4z_3)\\
x_\infty=\tau_0\tau_1((1+\mu^2s^4\omega^2)z_1+s^2(1+\omega^2)z_2
+s^4(1+\omega^4)z_3+s^6\omega^2(1+\omega^2)z_\infty)\end{array}\right.$$
\end{prop}
\dem The previous calculations show that one can form an $R$-basis
$\{Z_\bullet\}$ of $\Wr$ by setting
$$\begin{array}{lcl} Z_\infty=X_\infty &
\hspace{0.5cm} & Z_1=f\,X_\infty\\
Z_2=x_1x_2 X_\infty & & Z_3=(x_1+x_2)X_\infty\end{array}$$ The
expressions of the $X_\bullet$ in terms of the $Z_\bullet$ are
then  given by the equations (3.2) and (3.6). By duality, we
easily deduce the corresponding  change for coordinates.
$\square$

\section{Equations of $\widetilde{V}$ for a supersingular genus 2 curve in characteristic 2}

In [LP1] (section 5), it is shown that the morphism $D : M_X \to
\P^3=|2\Theta|$ is, as in the complex case (see [NR]), an
isomorphism when $X$ is an ordinary genus $2$ curve over an
algebraically closed field of characteristic $2$. As asserted in
[LP2], this identification extends to the relative case $\Xr \to
\Spec R$, so that the Frobenius morphism $\Xr \to \Xr_1$ induces,
by pull-back, a rational map

\begin{center}
\unitlength=0.6cm
\begin{picture}(14,5)
\put(3.4,3.5){$\P\, \Wr_1$} \put(10.3,3.5){$\P\, \Wr$}
\put(6.7,0.5){$\Spec R$} \put(5,3){\vector(1,-1){1.5}}
\put(10,3){\vector(-1,-1){1.5}} \put(6,3.7){\vector(1,0){3}}
\small{\put(7.5,3.9){$\Vr$}}
\end{picture}\end{center}
where $\Wr_1$ stands for the $2$-twist
$H^0({\JXr}_1,\,\O(2\Theta_1))$ of $\Wr$.

Over the generic point $\eta$, [LP1] (Proposition 3.1) gives the
form of the rational map \begin{eqnarray*} \Vr_\eta : \P\,
(\Wr_1)_\eta & \to & \P\,
(\Wr)_\eta\\
x=(x_\bullet) & \mapsto & (\lambda_BP_B(x) : \, \lambda_0P_0(x) :
\, \lambda_1P_1(x) : \, \lambda_\infty P_\infty(x))\end{eqnarray*}
where the $x_\bullet$ are the theta coordinates of the spaces
$(\Wr_1)_\eta$ and $(\Wr)_\eta$ (which correspond via the
$K$-semilinear isomorphism $i^* : (\Wr)_\eta \to (\Wr_1)_\eta$),
where the $P_\bullet$ are the quadrics
$$\begin{array}{lcl} P_B(x)=x_B^2+x_0^2+x_1^2+x_\infty^2
&\hspace{1cm} &
P_0(x)=x_Bx_0+x_1x_\infty\\
P_1(x)=x_Bx_1+x_0x_\infty & &
P_\infty(x)=x_Bx_\infty+x_0x_1.\end{array}$$ and where the
$\lambda_\bullet$ are nonzero constants depending on the curve
$X$. In [LP2] (section 3), we find  an explicit determination (up
to scalar) of these coefficients, namely
$$(\lambda_B:\,\lambda_0:\,\lambda_1:\,\lambda_\infty) =
(\sqrt{abc}:\,\sqrt{c}:\,\sqrt{b}:\,\sqrt{a})$$ Using the  formulae given in Proposition 3.7 and the expressions
(2.11), (2.12) and (2.13) for the coefficients $a$, $b$ and $c$, we can
compute the rational map
$$\Vr : z=(z_\bullet) \mapsto
(R_1(z):\,R_2(z):\,R_3(z):\,R_\infty(z))$$ First, we express the
polynomials $Q_\bullet(z)=P_\bullet(x)$ (notice that the
coefficients appearing in the formulae given in Proposition 3.7
are squared since we are dealing with elements of the $2$-twist
$(\Wr_1)_\eta$). We find the following:
$$ \begin{array}{rcl} \,Q_B(\underline{z})&=&
(\tau_0\tau_1)^4\left( (\nu_0\nu_1+\mu^4s^8\omega^4)^2z_1^2
+s^8(\nu_0\omega^4+\nu_1)^2z_2^2+s^{16}(\nu_0\omega^8+\nu_1)^2z_3^2+
s^{24}\omega^8(1+\omega^8)z_\infty^2\right)\\ \\

Q_0(\underline{z})&=&\tau_0^4\tau_1^2\left((\nu_0+\mu^2s^4\omega^2)^2z_1^2
+s^4\omega^4(\nu_0+\mu^2s^4)^2z_1z_2+s^8\omega^4(\nu_0\omega^2+\mu^2s^4)^2z_1z_3\right.\\
& & \hspace{0.5cm}
+s^8(1+\omega^4)z_2^2+s^{12}\omega^4(1+\omega^4)(z_2z_3+z_1z_\infty)+s^{16}(1+\omega^8)z^2_3\\
& & \hspace{1cm}
\left.+s^{16}\omega^4(1+\omega^4)(z_2z_\infty+s^4z_3z_\infty)\right)\\
\\

Q_1(\underline{z})&=&\tau_0^2\tau_1^4\left((\nu_1+\mu^2s^4\omega^2)^2z_1^2
+s^4(\nu_1+\mu^2s^4\omega^4)^2z_1z_2+s^8(\nu_1+\mu^2s^4\omega^6)^2z_1z_3\right.\\
& & \hspace{0.5cm}
+s^8\omega^4(1+\omega^4)z_2^2+s^{12}\omega^4(1+\omega^4)(z_2z_3+z_1z_\infty)+s^{16}\omega^8(1+\omega^8)z^2_3\\
& & \hspace{1cm}
\left.+s^{16}\omega^8(1+\omega^4)(z_2z_\infty+s^4\omega^4z_3z_\infty)\right)\\
Q_\infty(\underline{z})& =& (\tau_0\tau_1)^2s^8\omega^4
\left((\mu^4z_1^2
+z_2^2)+s^4(1+\omega^4)(z_2z_3+z_1z_\infty)+s^8\omega^4z^2_3\right)\\
\\
\end{array}$$

Secondly, using the coefficients $\lambda_\bullet$ and the
formulae of   Proposition 3.7 again (after inversion and up to
scalar), we obtain
 $$\begin{array}{rcl} R_1(z) & =&
\Frac{1}{(\tau_0\tau_1)^4s^{18}\omega^6(1+\omega)^6} \,Q_B(z)\\
R_2(z)&= &
\Frac{1}{(\tau_0\tau_1)^4s^2\omega^2(1+\omega^2)}\left[\Frac{1}{s^{18}\omega^6(1+\omega^2)}\,Q_B(z)
+\Frac{1}{s^8(1+\omega^4)}(Q_0(z)+Q_1(z))\right]\\
R_3(z)&= &
\Frac{1}{(\tau_0\tau_1)^4s^4\omega^2(1+\omega^2)}\left[\Frac{1}{s^{18}\omega^6(1+\omega^4)}\,Q_B(z)
+\Frac{1}{s^8\omega^2(1+\omega^4)}(\omega^2Q_0(z)+Q_1(z))\right]\\
R_\infty(z)& =&
\Frac{1}{(\tau_0\tau_1)^4s^6\omega^2(1+\omega^2)}\left[\Frac{1+\mu^2s^4\omega^2}{s^{18}\omega^6(1+\omega)^6}\,Q_B(z)\right.\\
& &\hspace{4cm} +\left.\Frac{1}{s^8\omega^4(1+\omega^4)}(\omega^4Q_0(z)+Q_1(z))+\Frac{1}{s^8\omega^4}Q_\infty(z)\right]\\
\end{array}$$
Now, using the expansions of $\nu_0$ and $\nu_1$, we obtain
\begin{eqnarray}\begin{array}{rcl} R_1(z) & =&
\mu^4z_1^2+z_2^2+s^2\bar{R}_1(z)\\
R_2(z)&= &
z_1^2+s^2\bar{R}_2(z)\\
R_3(z)&= &
(z_2z_3+z_1z_\infty)+s^2\bar{R}_3(z)\\
R_\infty(z)& =&
(\mu^4z_3^2+z_\infty^2+\mu^2z_1^2+z_1z_2)+s^2\bar{R}_\infty(z)\\
\end{array}\end{eqnarray}
where the $\bar{R}_\bullet(z)$ are quadrics with coefficients in
$R[\omega^{-1},\,(1+\omega)^{-1}]$.\\

\Rq A few more calculations show that the coefficients of the
quadrics $\bar{R}_\bullet(z)$ lie in fact in $R$, so that we can
choose any $\omega$ we want in $R\setminus\{0,\,1\}$ (and not necessarily
an element that specializes in $k\setminus\{0,\,1\}$). It means that any
deformation
of $X$ leads to the same equations.\\

Let us introduce the set of homogeneous coordinates
$\{y_1,\,y_2,\,y_3,\,y_\infty\}$ of $|2\Theta|$ given by
$$\begin{array}{lcl} y_1=z_1,&
\hspace{0.5cm} & y_2=z_2+\mu z_1,\\
y_3=z_3, & & y_\infty=z_\infty +\mu z_3.\end{array}$$ Thus, the
corresponding coordinates change in $|2\Theta_1|$ is obtained by
squaring the coefficients and we have proved the following.

\begin{theo} Let $X$ be a smooth, proper, and supersingular curve of genus $2$ over an  algebraically closed field of characteristic $2$. There exist
coordinates $\{z_\bullet\}$ (resp. $\{y_\bullet\}$) on $|2\Theta|$
(resp. $|2\Theta_1|$) such that the equations of $\widetilde{V}$
are given by
$$\widetilde{V} : |2\Theta_1| \to |2\Theta|,\hspace{0.5cm} y=(y_\bullet)
\mapsto z=(z_\bullet)=(Q_1(y):\,Q_2(y):\,Q_3(y):\,Q_\infty(y))$$
with
$$\begin{array}{lcl}
Q_1(y)=y^2_2, & \hspace{0.5cm} & Q_2(y)=y_1^2,\\
Q_3(y)=y_2y_3+y_1y_\infty, & & Q_\infty(y)=y_\infty^2+y_1y_2.\\
\end{array}\\ $$
\end{theo}

\section{Frobenius action on $M_X$ for a supersingular genus 2 curve in characteristic 2}

We know ([LP2], section 3) that the Kummer  surface
$\Kum_{\Xr_\eta}$ is defined, in terms of the theta
coordinates $\{x_\bullet\}$ on $|2\Theta|_\eta$, by the
homogeneous quartic
$$c(x_B^2x_0^2+x_1^2x_\infty^2)+b(x_B^2x_1^2+x_0^2x_²\infty^2)+a(x_B^2x_\infty^2+x_0^2x_1^2)+x_Bx_0x_1x_\infty,$$
 $a,\,b,\,c$ being the scalars appearing in the standard
birational model (2.3) of the curve $\Xr_\eta$.

The same kind of calculations as in section 4 give the equation of
$\Kum_X$ in the coordinate system $\{z_\bullet\}$ of
$|2\Theta|$, namely
\begin{eqnarray}
\mu^2z_1^3z_2+z_1^3z_\infty+z_1^2z_2z_3+\mu^4z_1^2z_3^2+z_1z_2^3+z_2^2z_\infty^2+z_3^4.
\end{eqnarray}
As in [LP1], we easily deduce from Theorem 4.2 and from the latter
calculations a complete description of the action of Frobenius on
$M_X$, more precisely of its separable part $V: M_{X_1}
\dashrightarrow M_X$, that we identify, using the isomorphism $D :
M_X \xrightarrow{\sim} |2\Theta|$ (see the introduction), with
$\widetilde{V} : |2\Theta_1| \dashrightarrow |2\Theta|$.

\begin{prop}Let $X$ be a smooth, proper, and supersingular curve of genus $2$ over an  algebraically closed field of characteristic $2$.
\begin{enumerate}
\item[1.] The semistable boundary of $M_X$ (resp. $M_{X_1}$) is
isomorphic (via $D$) to the Kummer quartic surface $\Kum_X$ (resp.
$\Kum_{X_1}$), an equation of which is (5.1). In particular, $V$
maps $\Kum_{X_1}$ onto $\Kum_X$.
 \item[2.] There is a unique stable bundle $E_{\rm bad}\in
M_{X_1}$ which is destabilized by Frobenius (i.e., $F^*E_{\rm
bad}$ is not semistable). We have $E_{\rm bad}=F_*B^{-1}$ and its
projective coordinates are $(0:0:1:0)$. \item[3.] Let $H_1$ be the
hyperplane in $|2\Theta_1|$ defined by $y_2=0$. The map $V$
contracts $H_1$ to the conic $\Kum_X\cap H$, where $H$ is the
hyperplane in $|2\Theta|$ defined by $z_1=0$.\item[4.] The fiber
of $V$ over a point $[E]\in M_X$ is
\begin{enumerate}
 \item[$\bullet$] a single point $[E_1]\in
M_{X_1}$, if $[E]\notin H$,
\item[$\bullet$] empty, if $[E]\in H
\setminus \Kum_X \cap H$,
\item[$\bullet$] a projective line
passing through $E_{\rm bad}$, if $[E]\in \Kum_X \cap
H$.\end{enumerate}
 In particular, $V$ is dominant, nonsurjective.
 \item[5.]  The total (resp. separable) degree of
$V$ is $4$ (resp. $1$).
\end{enumerate}
\end{prop}
\dem The fact that the Kummer  quartic surface $\Kum_X$ is the
semistable boundary of $M_X$ comes from [NR]. On the one hand, the
inverse image of $\Kum_X$ is a closed subspace in $|2\Theta_1|$
the ideal of which is generated by the homogeneous polynomial
obtained by replacing the $\{z_\bullet\}$ in equation (5.1) by the
$\{Q_\bullet(y)\}$ given in Theorem 4.2, namely
$$y_2^4[y_1^3y_\infty+\mu^2y_1^2y_2^2+y_1^2y_2y_3+\mu^4(y_1^2y_\infty^2+y_2^2y_3^2)+y_1y_2^3+y_2^2y_\infty^2+y_3^4]$$
On the other hand, we obtain the equation
$$y_1^3y_\infty+\mu^2y_1^2y_2^2+y_1^2y_2y_3+\mu^4(y_1^2y_\infty^2+y_2^2y_3^2)+y_1y_2^3+y_2^2y_\infty^2+y_3^4$$
of $\Kum_{X_1}$ in $|2\Theta_1|$ by squaring the
coefficients in equation (5.1) and replacing the coordinates
$\{z_\bullet\}$ on $|2\Theta_1|$ by their expressions
$$\begin{array}{lcl} z_1=y_1&
\hspace{0.5cm} & z_2=y_2+\mu^2y_1\\
z_3=y_3 & & z_\infty=y_\infty +\mu^2y_3\end{array}$$ in terms of
the $\{y_\bullet\}$. Thus, $V^{-1}(\Kum_X)=\Kum_{X_1}\bigcup H_1$
(where $H_1=(y_2=0)$), hence (1). Furthermore, $H_1$ is mapped
onto the hyperplane $H=(z_1=0)$ of $|2\Theta|$, hence (3).

In [LP1] (Proposition 6.1), it is shown that there is exactly one
base-point, namely $F_*B^{-1}$. Now, solving $Q_\bullet(y)=0$, we
find that $(0:\,0:\,1:\,0)$ is the unique base-point of $V$, hence
(2).

Let $[E]$ be a $k$-point of $M_X$, with coordinates $(a^2:\,b^2:\,
c^2:\,d^2)$. Let us solve the system
\begin{eqnarray}
y_2^2=a^2,\hspace{5mm}
y_1^2=b^2,\hspace{5mm}y_2y_3+y_1y_\infty=c^2,\hspace{5mm}y_\infty^2+y_1y_2=d^2.
\end{eqnarray}
We must have $y_1=b$ and $y_2=a$, so $y_\infty=d+\sqrt{ab}$ and
finally $ay_3=c^2+bd+b\sqrt{ab}$. If $[E]$ belongs to $H$, i.e.,
if $a=0$, it has solutions if and only if $c^2=bd$, that is if and
only if $[E]$ belongs to $H \cap \Kum_X$. In that case,
there is no condition on $y_3$ and the inverse image of such a
point is a projective line passing through $E_{\rm bad}$. If $[E]$
does not belong to $H$, i.e., if
$a\neq 0$), we see that the inverse image is a single point.

Finally, defining $u=y_1/y_2$, $v=y_3/y_2$ and $w=y_\infty/y_2$,
the field extension $$k(u^2,\,v+uw,\,w^2+u) \subseteq k(u,\,v,\,w)
\cong k(u^2,\,v+uw,\,w^2+u)[t,\,s]/(t^2-u^2,\,s^2-w^2)$$
corresponds to the Verschiebung $V$, and the last
assertion follows. $\square$\\

Denote by $N_X$ (resp. $N_{X_1}$) the moduli of semistable
bunbles with rank $2$ and degree $0$ over $X$ (resp. $X_1$). In the
case of an ordinary genus $2$ curve, [LP1] (Proposition 6.4) showed
the surjectivity of the (rational map) Verschiebung $N_{X_1}
\dashrightarrow N_X$, $ [E_1] \mapsto [F^*E_1]$. In the case of a
supersingular curve, this result does not hold any more. Let $E$
be a vector bundle in $N_X$. The Jacobian $JX$ being divisible, we
can assume that $\det E$ is trivial. Now, if
$[F^*E_1]=[E]$, we have $F^*(\det E_1)=\det E$ and
$\det E_1$ has to be a $2$-torsion point of $JX_1$.
Therefore, $\det E_1$ is trivial as well and the following
is a corollary of Proposition 5.1.

\begin{prop} Let $X$ be a smooth, proper, and supersingular curve of genus $2$ over an  algebraically closed field of characteristic $2$. The rational map
$N_{X_1} \to N_X$ given by $[E_1] \mapsto [F^*E_1]$ is not
surjective.
\end{prop}

We are now interested in those vector bundles over $X$ that are
destabilized by a finite number of iteration of the (absolute)
Frobenius. Let $\Omega^{\rm Frob}$ be the complementary set of
classes of semistable rank $2$ vector bundles $E$ with trivial
determinant over $X$ such that $F_{\rm abs}^{(n)\,*}E$ is
semistable for all $n\geq 1$.

\begin{prop}Let $X$ be a smooth, proper, and supersingular curve of genus $2$ over an  algebraically closed field of characteristic $2$. The open subset
$M_X\setminus\{E_{\rm bad}\}$ is stable under the action of the (absolute)
Frobenius. In particular, $\Omega^{\rm Frob}$ is the Zariski open
(dense) subset $M_X\setminus \{E_{\rm bad}\}$ of $M_X$.
\end{prop}
\dem Pulling back a semistable bundle $E$ over $X$ by the
semilinear isomorphism $i : X_1 \to X$, we obtain a semistable
bundle $i^*E$ over $X_1$. If $E$ has trivial determinant, $i^*E$
has trivial determinant as well. If $E$ ($\neq E_{\rm bad}$) has
coordinates $(a:\, b:\, c:\,d)$ in the system $\{z_\bullet\}$,
$i^*E$ has squared coordinates in the corresponding system
$\{z_\bullet\}$ of $M_{X_1}$. Using Theorem  4.2  (actually, the
quadrics given by equations (4.1) from which we deduce Theorem
4.2 after the indicated coordinates change), we find that
$F_{\rm abs}^*E=V^*(i^*E)$ has coordinates
$$(\mu^4a^4+b^4:\,a^4:\,a^2d^2+b^2c^2:\,\mu^4c^4+d^4+\mu^2a^4+a^2b^2)$$
in the system $\{z_\bullet\}$, hence cannot be $E_{\rm bad}$. $\square$\\

\begin{cor} Let $k$ be  a finite field of characteristic $2$, let
 $X$ be a smooth, proper, and supersingular curve of genus $2$ over   $k$, and let $E$ be a semistable rank $2$
vector bundle with trivial determinant, defined over $k$ and
different from $E_{\rm bad}$. Then some twist ${F^{(n_0)}}^*E$
comes from a continuous representation of the algebraic
fundamental group $\pi_1(X)$ in $\GL_2(\bar{k})$.\end{cor}
\dem For such a vector bundle, the sequence $F_{\rm abs}^{(n)\,*}E$
($n\geq 1$) takes its values in the finite set $M_X(k)$ of $S$-equivalence classes of
$k$-rational rank $2$ semistable vector bundles with trivial
determinant. Therefore, the sequence
$F_{\rm abs}^{(n+n_0)\,*}E$ ($n\geq 1$) is periodic for suitable $n_0$
and we use [LS] (Theorem 1.4)  to conclude. $\square$

\vspace{1cm}

\noindent \Large \textbf{References}\\
\normalsize
\begin{enumerate}
\item[{[AG]}] \textsc{J. Arledge}, \textsc{D. Grant}: \emph{An
explicit theorem of the square for hyperelliptic Jacobians},
Michigan Math. J.  {\bf49} (2001), 485--492.

\item[{[DO]}] \textsc{I. Dolgachev}, \textsc{D. Ortland}:
\emph{Point sets in projective space and theta functions},
Ast\'erisque {\bf169}, 1988.

\item[{[H]}] \textsc{R. Hartshorne}: \emph{Algebraic geometry},
Graduate Texts in Mathematics {\bf52}, Springer, New-York, 1977.

\item[{[L]}] \textsc{H. Lange}: \emph{\"{U}ber die Modulschemata
der Kurven vom Geschlecht $2$ mit $1$, $2$ oder $3$ Weierstrasspunkten},
J. Reine Angew. Math. {\bf277} (1975), 27--36.

\item[{[LP1]}] \textsc{Y. Laszlo}, \textsc{C. Pauly}: \emph{The
action of the Frobenius map on rank 2 vector bunbles in
characteristic 2}, J. of Alg. Geom. {\bf11} (2002), 219--243.

\item[{[LP2]}] \textsc{Y. Laszlo}, \textsc{C. Pauly}: \emph{The
Frobenius map, rank 2 vector bunbles and Kummer's quartic surface
in characteristic 2 and 3}, Advances in Mathematics {\bf185}
(2004), 246--269.

\item[{[LS]}] \textsc{H. Lange}, \textsc{U. Stuhler}:
\emph{Vektorb\"{u}ndel auf Kurven und Darstellungen der
algebraischen Fundamentalgruppe}, Math. Zeit. {\bf156} (1977), 73--83.

\item[{[M1]}] \textsc{D. Mumford}:  \emph{Abelian varieties},
Tata Institute of Fundamental Research Studies in Mathematics {\bf
5}, Bombay, 1970

\item[{[M2]}] \textsc{D. Mumford}:  \emph{On equations defining
abelian varieties. I.}, Invent. Math. {\bf1} (1966), 287--354.

\item[{[NR]}] \textsc{M.S. Narasimhan, S. Ramanan}: \emph{Moduli
of vector bundles on a compact Riemann surface}, Ann. of Math. {\bf89}
 (1969), 14--51.

\item[{[R]}] \textsc{M. Raynaud}: \emph{Sections des fibr\'es
vectoriels sur une courbe}, Bull. Soc. Math. France {\bf110}
(1982), 103--125.

\item[{[S]}] \textsc{T. Sekiguchi}: \emph{On projective
normality of abelian varieties. II.}, J. Math. Soc. Japan {\bf29}
(1977), 709--727.

\end{enumerate}

\vspace{2cm}

\noindent Laurent Ducrohet\\
Universit\'e Pierre et Marie Curie\\
Analyse alg\'ebrique, UMR 7586\\
4, place Jussieu\\
75252 Paris Cedex 05 France\\
e-mail: ducrohet@math.jussieu.fr

\end{document}